\documentclass[aps,preprint,a4paper,nofootinbib,floatfix]{revtex4}

\usepackage{graphicx,color}
\usepackage{amsmath,amssymb,amsfonts,amsthm,amscd,bm}
\usepackage{booktabs}
\usepackage{cancel}

\usepackage{amsmath}
\usepackage{amssymb}
\usepackage{amsfonts}
\usepackage{amsthm}
\usepackage{amscd}
\usepackage{bm}
\usepackage{graphicx,color}
\usepackage{booktabs}
\usepackage{listings}
\usepackage{xcolor,colortbl}
\usepackage{float}
\usepackage{mathtools}

\def\bar{\overline}

\def\*{\star}
\def\[{\left[}
\def\]{\right]}
\def\({\left(}      

\def\){\right)}

\def\frac#1#2{\dfrac{#1}{#2}}
\def\inv#1{\dfrac{1}{#1}}

\def\2pi{\hbox{$2\pi i$}}

\def\dsl{\raise.15ex\hbox{/}\kern-.57em\partial}
\def\Dsl{\,\raise.15ex\hbox{/}\mkern-.13.5mu D}
   \def\CN{{\cal N}}   
\def\CP{{\cal P}}

\def\2pi{\hbox{$2\pi i$}}

\def\dsl{\raise.15ex\hbox{/}\kern-.57em\partial}
\def\Dsl{\,\raise.15ex\hbox{/}\mkern-.13.5mu D}
\font\numbers=cmss12
\font\upright=cmu10 scaled\magstep1
\def\stroke{\vrule height8pt width0.4pt depth-0.1pt}
\def\topfleck{\vrule height8pt width0.5pt depth-5.9pt}
\def\botfleck{\vrule height2pt width0.5pt depth0.1pt}
\def\Zmath{\vcenter{\hbox{\numbers\rlap{\rlap{Z}\kern
    0.8pt\topfleck}\kern 2.2pt
    \rlap Z\kern 6pt\botfleck\kern 1pt}}}
\def\Qmath{
    \vcenter{\hbox{\upright\rlap{\rlap{Q}\kern3.8pt\stroke}\phantom{Q}}}}
\def\Nmath{\vcenter{\hbox{\upright\rlap{I}\kern 1.7pt N}}}
\def\Cmath{\vcenter{\hbox{\upright\rlap{\rlap{C}\kern
                   3.8pt\stroke}\phantom{C}}}}
\def\Rmath{\vcenter{\hbox{\upright\rlap{I}\kern 1.7pt R}}}
\def\Z{\ifmmode\Zmath\else$\Zmath$\fi}
\def\Q{\ifmmode\Qmath\else$\Qmath$\fi}
\def\N{\ifmmode\Nmath\else$\Nmath$\fi}
\def\C{\ifmmode\Cmath\else$\Cmath$\fi}
\def\R{\ifmmode\Rmath\else$\Rmath$\fi}
\def\barray{\begin{eqnarray}}
\def\earray{\end{eqnarray}}
\def\beq{\begin{equation}}
\def\eeq{\end{equation}}

\def\AA{\leavevmode\setbox0=\hbox{h}
\dimen0=\ht0 \advance\dimen0 by-1ex\rlap{\raise.67\dimen0\hbox{\char'27}}A}

\def\Arg{{\rm Arg}\,}

\def\prob{{ \bf Pr}}
\def\Ex{{\bf E}}

\def\Primes{\mathbb{P}}
\def\Primesp{\mathbb{P}'}
\def\Pensemble{{\bf \CP}}

\def\dist{{\overset{d}{\longrightarrow}}}

\makeatletter
\def\iddots{\mathinner{\mkern1mu\raise\p@
\vbox{\kern7\p@\hbox{.}}\mkern2mu
\raise4\p@\hbox{.}\mkern2mu\raise7\p@\hbox{.}\mkern1mu}}
\makeatother

\theoremstyle{plain}
\newtheorem{theorem}{Theorem}

\newtheorem{corollary}{Corollary}

\theoremstyle{remark}
\newtheorem{remark}{Remark}

\begin{document}

\title{
 Central Limit Theorems  for series of Dirichlet characters
}
\author{
 Andr\'e  LeClair\footnote{andre.leclair@gmail.com}
}
\affiliation{Cornell University, Physics Department, Ithaca, NY 14850} 

\begin{abstract}

For a given   Dirichlet character   $\chi (n) = e^{i \theta_n}$,  we prove central limit theorems for the series 
$\sum_{p'}  \cos \theta_{p'}$ for non-principal characters,  and $\sum_{p' }  \cos (t \log p')$  for principal characters, 
where $p'$  are integers based on a variant of Cram\'er's  random model for the primes.     For non-principal characters,  we use these results to show 
that the Generalized Riemann Hypothesis  for the associated $L$-function  is true with probability equal to one.     For principal characters 
we propose  how to extend these arguments to $\Re (s)  =  t \to \infty$.  

\bigskip\bigskip
{\it 
In memory of my daughter Alexandra LeClair,   who passed away during the course of this work.
}

\end{abstract}

\maketitle

\section{Introduction}

This article concerns the growth of certain infinite series defined over prime numbers.   
In this Introduction  we explain the motivation for this study,   define the series in question,  and briefly summarize
some  of our results.   

  Let $\chi (n)$  denote a Dirichlet character,   and as usual   
 define the Dirichlet series 
\beq
\label{Dir}  
L(s, \chi ) = \sum_{n=1}^\infty    \frac{\chi (n)}{n^s} 
\eeq
where $s= \sigma +i t$ is a complex variable.    The Riemann zeta function $\zeta (s)$  corresponds to the $L$-function 
for the trivial character modulo $1$ where all $\chi (n) =1$.     Due to the completely multiplicative property of the characters,
$L$  enjoys an Euler product formula:
\beq
\label{EPF}  
L(s, \chi)   =  \prod_{n=1}^\infty  \( 1 -  \frac{\chi (p_n)}{p_n^s} \)^{-1} 
\eeq
where $p_n$ is the $n$-th prime.       The above formula is known to be valid for $\Re (s) >1$  where both sides of the equation 
converge absolutely.    Using the  validity of the Euler product one can easily see that there are no zeros with $\Re (s) > 1$;   in particular,  
$\log L$ is finite in this region since the series converges.       If somehow the Euler product were valid for $\Re (s) > 1/2$,   then the Generalized Riemann Hypothesis would
follow by the same argument together with the functional equation that relates $L(s,\chi)$ to $L(1-s, \bar{\chi})$.     For Riemann $\zeta$ itself,  and other $L$-functions
based on {\it principal}  Dirichlet characters,   it is well understood that the Euler product in the  above form  is  not  valid for $1/2 <  \Re (s) \leq 1$,  essentially due to the pole
at $s=1$:    domains of convergence of 
Dirichlet series are always half-planes,   and  due to this pole,  
the logarithm must be divergent for $\Re (s) \leq 1$.      However for {\it non-principal} characters,  the associated $L$-functions have no pole at $s=1$,   and thus
the validity of the Euler product for $\Re (s) >1/2$  is  theoretically possible,  although difficult to study.     

For the moment,  let us not distinguish between the cases of principal verses non-principal characters,  although for the reasons discussed above there will subsequently 
be significant differences.  
Let the character have modulus $k$.      The character $|\chi (n)| = 1$ if  $(n,k)=1$,  i.e. $n,k$ are coprime,  otherwise $\chi (n) = 0$.   
The non-zero characters are all roots of unity,   
so let us define the angles $\theta_n$: 
\beq
\label{thetan}
\chi (n) =  e^{i \theta_n }, ~~~~~\forall ~\chi(n) \neq 0
\eeq   
Now consider 
the series 
\beq
\label{BNt} 
B_N (t, \chi)  =  \sum_{n=1}^N   \cos \( \theta_{p_n}  +  t \log p_n  \)  
\eeq
where 
it is implicit that the terms corresponding to the finite number of primes for which  $\chi(p_n)=0$ are omitted in the sum.  
In  \cite{EPF1,EPFchi} it was proven that if 
 $B_N = O(\sqrt{N})$  as  $N \to \infty$,   then the Euler product formula is valid for $\Re (s) > 1/2$ because it converges in this region.     
The proof involved an Abel transform for the logarithm of
the Euler product,  the Prime Number Theorem (PNT),  and a bound on the sum over gaps between primes.     Let us now specialize to non-principal characters. 
Because of the half-plane convergence property mentioned above,  to establish validity of the Euler product formula for $\Re (s) > \sigma_c$ for some $\sigma_c$,   
it is sufficient to prove convergence at a single value of $t$.    Since there is no pole
at $s=1$,  the simplest  choice is $t=0$.    It is sufficient then to consider the series 
\beq
\label{CNdef} 
C_N =  \sum_{n=1}^N   \cos   \theta_{p_n}
\eeq

As stated above,   a proof that $C_N = O(\sqrt{N})$  would establish the validity of the Generalized Riemann Hypothesis for all non-principal characters.  
It is a completely deterministic series which depends on the actual primes which are largely unknown for large $N$,   thus it is difficult,  if not impossible,  to 
compute it for large enough $N$.     However if one is only interested in its growth as a function of $N$,  since $C_N$ is a series,  the fluctuations coming from 
the precise values of individual  primes may not be important for  determining this growth.   In other words,   the growth of $C_N$  may only depend on some 
global properties of the set of primes,  such as their average spacing,  etc.           In \cite{EPFchi}  it was conjectured that $C_N = O(\sqrt{N})$  
based on the heuristic argument that it behaves like a random walk due to the multiplicative independence  of the primes.  
     In this article we apply methods of probability theory to further study this problem.     The idea of using probability methods in number theory 
     is certainly not new,  and at least goes back to work of Cram\'er \cite{Cramer},  which we will utilize.          

Let $\Primes = \{p_1, p_2, \ldots \}$  denote the  set of primes,  where $p_1 =2,  p_2 = 3, $ and so forth. 
     We will consider replacing $\Primes$ with the set $\Primesp = \{ p'_1, p'_2, \ldots \}$,  
which are a random,  independent, ordered sequence of  integers,   and  will study  $C'_N =  \sum_{n=1}^N   \cos \theta_{p'_n}$.   
 The $p'_n$  will  be  constrained  to  satisfy some global properties of the known primes,   to be  specified below. 
    Since the $p'_n$ are now random variables,
we will  consider $\Pensemble =  \{ \Primesp \}$   which is the ensemble of all possible  $\Primesp$,  i.e the set of sets $\Primesp$.      We will refer to
$\Pensemble$  as the pseudo-prime ensemble,   and a specific element $\Primesp \in \Pensemble$ as a state of this  ensemble.   The actual primes 
$\Primes$ are then simply one state in the pseudo-prime ensemble.      This terminology  is borrowed from statistical mechanics in physics.    For instance,  
for a gas of free particles,   the states correspond to specific values for the positions and velocities of every individual particle,   and  various canonical ensembles are the
set of all such states subject to certain constraints,  such as the total number of particles or total energy held fixed.     In light of this analogy,   the idea we are pursuing here
is that,  for example,  the macroscopic pressure of a gas of a large number of particles hardly depends on what specific state they are in, which is unknowable,  
 and for the same reasons  we expect the 
global (macroscopic)  properties of  $C'_N$,  in particular its growth as a function of $N$,   does not depend on the detailed properties of $\Primesp$.     The aim of this article to 
make such statements precise using the theory of probability,   as in statistical mechanics.    

We  need to be specific about the ensemble $\Pensemble$ and its probability measure.     We will require two   properties.    The primes $p_n$ are independent, 
more specifically they are multiplicatively independent.    We thus require the $p'_n$  to also  be independent.     Secondly,  we require that the counting of $p'$ is essentially equivalent 
to that implied by the Prime Number Theorem.     Namely, as usual  let $\pi (x)$ be the number of primes less than $x$.    The PNT gives the 
leading behavior  $\pi (x)  \approx  {\rm Li} (x) \approx  x / \log x$.     Let $\pi' (x)$ be the analogous quantity for the primes $p'_n$,  i.e. the number of $p' <x$.     This counting function is now a random
variable,  and we require that its expectation value is approximately the leading term in $\pi (x)$:
\beq
\label{Epi} 
\Ex [ \pi'(x) ]  \approx \pi(x) \approx  {\rm Li} (x) \approx   \frac{x}{\log x} 
\eeq

There are many possible choices of $\Pensemble$ compatible with these requirements.   One particular interesting one is to 
take $p'_n$  to be a random integer satisfying  $p_n \leq p'_n \leq p_{n+1}$.    Remarkably,   Grosswald and Schnitzer \cite{Grosswald}  proved that
if one defines a $\zeta$ function via an Euler product from the  $\{ p'_n \}$,  as in \eqref{EPF},    then all of these possible  $\zeta$'s can be analytically continued into the 
critical strip and  have the {\it same}  zeros as
Riemann $\zeta$ there.   They proved a similar result for Dirichlet $L$-functions which will be discussed below.    
   For our purposes however,  this choice is more difficult to analyze than necessary.    Instead we will use a variant  of the  Cram\'er model \cite{Cramer}
   which depends on the modulus of the character $\chi$.    For instance, in the simplest case of modulus $k=1$,     
 for each integer $n$,  the probability that $n\in \Primesp$ equals $1/\log n$.     
We will then prove that  $C'_N$  obeys a central limit theorem,  i.e. when properly normalized,  it has a normal distribution:

\medskip
\begin{theorem}\label{CLTC}   For non-principal characters of modulus $k$, 
\beq 
\label{CLTforC}
\sqrt{\frac{1 +  \tfrac{\log \log N}{\log N} }{s^2 N}}  ~C'_N  ~~ \dist ~~ \CN (0,1) 
\eeq
with  
\beq
\label{s2} 
s^2 =  a \, \frac{\varphi (k)}{k}
\eeq
where $a=1$ if the characters $\chi$ are all real, (i.e. all $\pm1$),  otherwise $a=1/2$,  $\varphi (k)$ is the Euler totient,  and $\CN (\mu, \sigma )$ is the normal distribution
with  mean  $\mu$ and standard deviation $\sigma$.    
 \end{theorem}  
 
 As $N\to \infty$,  the $\tfrac{\log \log N} {\log N}$ can of course be neglected,  however we retain  it in order to provide numerical evidence 
 at large but finite $N$.    We will use this theorem to say something precise about the growth of the original series $C_N$.

For principal characters all the angles $\theta_{p_n}$ are zero and one needs to consider now the series 
\beq
\label{BNdef} 
B_N (t)   = \sum_{ {n=1} \atop  {(p_n ,k)=1}}^N  \cos (t \log p_n)
\eeq
Here,  obviously we are interested in  $t\neq 0$,  which as explained above, in relation to the validity of the Euler product,   this  is due to   the pole in $\zeta (s)$ at $s=1$.     
As before we define $B'_N (t)$  as above with $p_n \to p'_n$.   Below we will prove a central limit theorem   for $B'_N (t)$  in the limit of 
large $t$  (Theorem \ref{CLTB} below).

\section{Non-principal case}   

As explained in the Introduction,   given a non-principal  Dirichlet character $\chi$ of modulus $k$,   we are interested in the series 
\beq
\label{CNp}
C'_N  =  \sum_{ {n=1} \atop {(p'_n ,k) =1} }^{N}    ~ \cos \theta_{p'_n}
\eeq
where the angles $\theta_n  = \Arg \chi (n)$,  and   $\{ p'_n \} = \Primesp$ is one state in  the ensemble $\Pensemble$ appropriate to Cram\'er's model.   

We first describe how to implement the Cram\'er model and  generate the states $\Primesp$ in a way that is 
simple to study both analytically and numerically.     For simplicity we exclude $p'_1 =  2$ from $\Primesp$ .
This does not affect the large $N$ result we will obtain.        
For each  $n\geq 3$,   let $r_n$ be a random variable uniformly distributed on the interval $[0,1]$,  and define $z_n  $ as follows
\beq
\label{zndef}
z_n  = 
\begin{cases*}
1 &\text{ if}   $r_n \leq \tfrac{1}{\log n} $   \\
  0  &\text{otherwise} 
  \end{cases*}
\eeq
Then,  by definition, for $n\geq 3$,  $n\in \Primesp$  if $z_n = 1$.   
The $z_n  $ are independent  random variables,   with probabilities  $\prob[z_n =1] = 1/\log n$.   
 We have excluded $n=2$ since $1/\log 2 > 1$.    
 The counting formula $\pi' (x)$  for the number of $p'  \leq x$ 
is then simply 
\beq
\label{pip}
\pi' (x)  =  \sum_{n \leq x}    z_n  
\eeq
Since $\Ex[z_n] = 1/\log n$, 
\beq
\label{Expip}
\Ex[\pi' (x) ]   =   \sum_{n=3}^x   \inv{\log n}   \approx  \int_3^x   \frac{du}{\log u}    \approx \frac{x}{\log x} 
\eeq
in accordance with \eqref{Epi}.   

\def\nN{n_N} 
\def\znk{z_{n,k}}

In order to implement $(p'_n , k) =1$ in \eqref{CNp},  let us slightly modify the definition \eqref{zndef} to the following 
\beq
\label{znkdef}
z_{n,k}    = 
\begin{cases*}
1 &\text{ if}   $r_n \leq \tfrac{1}{\log n} $   ~ {\rm and}   ~ (n,k) =1\\
  0  &\text{otherwise} 
  \end{cases*}
\eeq
for $n\geq 3$.     
The series $C_N'$ defined in the Introduction is now modeled as 
\beq
\label{CNz}
C'_N  =  \sum_{n=1}^{p'_N}    \znk  \cos \theta_n 
\eeq
Summing over all integers up to $p'_N$ ensures there are $N$ terms in the sum when $k=1$.  
The series $C'_N$ is now a random variable with a well-defined probability distribution.    Let us now prove Theorem  \ref{CLTC}.

\begin{proof}    
\label{proofCLTC} 
{\it (of Theorem \ref{CLTC})}.   

Let us write  $C'_N =  \sum_{n}  c_n$  where $c_n =  \znk  \cos \theta_n$.      The $c_n$ are independent random variables
however they are not identically distributed,  and thus the classical (Lindeberg-L\'evy)  central limit theorem (CLT) does not apply.     
However  Lyapunov's CLT does.     More generally,  let $x_n$, $n=1,2,\ldots, N$  be independent random variables with finite  mean $\mu_n$  and variance $\sigma_n^2$, 
 which are allowed to vary with $n$,  and define  the series $X_N = \sum_{n=1}^N  x_n$.       
Define  $m_N$ as  the expectation value  of $X_N$,
\beq
\label{mN}
m_N  =  \Ex \Bigl[X_N \Bigr]  =  \sum_{n=1}^{N}  \mu_n, 
\eeq 
and  $s_N^2$ the sum of variances
\beq
\label{sN2} 
s_N^2  =  \sum_{n=1}^{N}  \sigma_n^2 
\eeq
 If the Lyapunov condition is satisfied,  namely if for some $\delta >0$  
 \beq
 \label{LyCond}  
 \lim_{N \to \infty}    \inv{s_N^{2 + \delta} }  \sum_{i = 1}^{p'_N}   \Ex \Bigl[ | x_n - \mu_n |^{2 + \delta} \Bigr]   = 0,  
 \eeq 
 then Lyapunov's theorem states that 
\beq
\label{Lyapunov}  
\inv{s_N}     \Bigl(  X_N  - m_N \Bigr)   ~ \dist ~  \CN(0,1)
\eeq

Let us now apply this to $X_N=C'_N$.  
First consider $m_N$.   
For non-principal characters,   one has 
\beq
\label{chisum}
\sum_{n=1}^{k-1}   \chi(n) =0 
\eeq
Thus the angles $\theta_n$ are equally spaced on the unit circle.   If the $p'_n$ are random,   
then  $\sum \cos \theta_{p'_n}$ is  always close to zero,  and on average is zero.        
 We will only need the weaker statement that  $m_N =O(1)$.    

Let us now turn to $s_N^2$:
\beq
\label{sN2}
s_N^2  =  \sum_{n=3}^{p'_N}   \Ex \[   \znk^2  \cos^2 \theta_n  \]
\eeq
Let us invoke an Abel transform  (integration by parts)  to re-express the above series in terms of $\sum \cos^2 \theta_n$.     
If the characters are all real,  then $\cos^2 \theta_n =1$  for all $n$.   On the other hand if they are complex and equally spaced on the unit circle
the average of $\cos^2 \theta_n$ is $1/2$.    Let us distinguish these two cases by  defining $a=1$ and $a=1/2$ respectively.   
For  $1\leq n \leq k$,  there are exactly   $\varphi (k)$ non-zero characters $\chi (n)$.    Since the characters are periodic,  
$\chi (n + k) = \chi (n)$,   the fraction of non-zero terms in the above sum is $\varphi (k)/k$.    
One clearly has   $\Pr \[ \znk^2 = 1 \] = 1/\log n$,  which implies  
\beq
\label{sN2b}
s_N^2  =      s^2   \sum_{n=3}^{p'_N}   \inv{\log n}       
\eeq
where $s^2$ is defined in \eqref{s2}.   
Next we use $p'_N \approx N \log N$ to obtain in the limit of large $N$:
\beq
\label{sN2c}
s_N^2  \approx  
s^2   \int_3^{N\log N}    \frac{du}{\log u}   \approx  s^2  N  \(1+ \frac{\log \log N}{\log N} \)^{-1}
\eeq
The Lyapunov condition is  easily verified for integer $\delta$,   since $s_N = O(\sqrt{N})$ and the expectation in \eqref{LyCond}  is $O(N)$ for any $\delta$.         
The theorem then follows from Lyapunov's result \eqref{Lyapunov},   using $\lim_{N\to \infty}  m_N/s_N = 0$.

\end{proof}

\begin{figure}[t]
\centering\includegraphics[width=.5\textwidth]{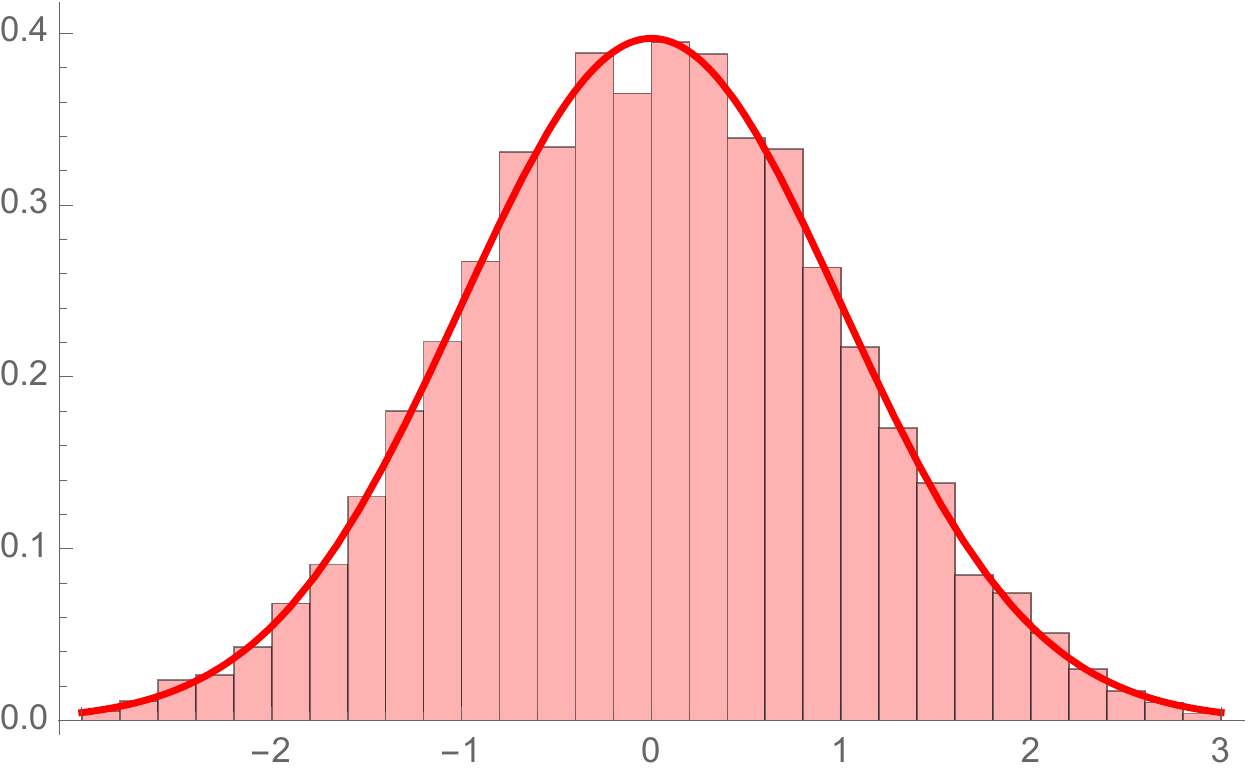}
\caption{Numerical evidence for Theorem \ref{CLTC}  based on the character \eqref{chi7}.  We fixed N=5,000.  
Displayed is a normalized histogram for $10,000$ states $\Primesp$.     The red curve is the  fit to the data,  which is the normal 
distribution $\CN(0.000500253, 1.0051)$.   The  blue curve is the prediction $\CN(0,1)$ which is nearly invisible since
it is indistinguishable from the fit.  
  }
\label{CLTCfig}
\end{figure}

In Figure \ref{CLTCfig}  we present compelling numerical evidence for Theorem \ref{CLTC}.  
We chose the following character   with $k=7$:
\beq
\label{chi7}
 \chi(1), \ldots, \chi(7) =   
1, e^{ 2\pi i /3},e^{\pi i / 3}, e^{-2\pi i / 3}, e^{-\pi i / 3}, -1,0
\eeq
Here,  $a=1/2$ and $\varphi(7)=6$.   We fixed $N=5,000$ and  generated  $10,000$  states $\Primesp$  numerically according to \eqref{znkdef};  displayed is a normalized histogram.     Performing a fit to a normal distribution 
gave $\CN(0.000500253, 1.0051)$. 

Although the true primes are obviously special,   they fall well within the bell curve,  which is to say they  are rather  ``normal".       Namely,  
for the special state  $\Primesp$  equal to the actual primes $\Primes$,   the  LHS of 
\eqref{CLTforC}  for the series $C_N$  equals  $-0.145$ for $N=5000$.    

 \bigskip
 
 \begin{theorem}   
 \label{CNsqrt}
  For any $\epsilon >0$,  in the limit $N \to \infty$,    $$C'_N = O(N^{1/2+ \epsilon})$$  with probability equal to $1$.  
 
 \bigskip
 \begin{proof}  
 Using the normal distribution of Theorem \ref{CLTC},  
 in the limit $N \to \infty$ one has  
 \beq
 \label{PrCN}
 \prob \[  C'_N   \leq  s \kappa N^{1/2+ \epsilon}   \]   =  
 1 -  \frac{e^{- \kappa^2 N^{2\epsilon}}} {\sqrt{2 \pi}  \kappa N^\epsilon}  \(  1 - O\(  \inv{\kappa^2  N^{2\epsilon} } \)  \)
 \eeq
 For any $\epsilon > 0 $,   
 \beq
 \label{ProbOne}
 \lim_{N \to \infty}  \prob \[  C'_N = O(N^{1/2+ \epsilon}) \]  = 1    
 \eeq
 
\end{proof}
 \end{theorem}

Given any particular state $\Primesp$,  we can define  the function 
\beq
\label{Lprime} 
L' (s, \chi)  =  \prod_{n=1}^\infty  \(  1 - \frac{\chi (p'_n)}{({p'_n})^s} \)^{-1}
\eeq

\begin{theorem}
\label{RH1}  
With probability equal to $1$,  all the functions $L'$ have no zeros with $\Re (s) > 1/2 + \epsilon$ for any $\epsilon > 0$.    

\begin{proof}  
Consider the limit  $\epsilon \to 0^+$ in Theorem \ref{CNsqrt}.    
It was shown in \cite{EPFchi}  that if $C'_N = O(N^{1/2+\epsilon})$,   then the logarithm of the product on the RHS of \eqref{Lprime} converges for $\Re (s) > 1/2+\epsilon$.  
Thus the very definition of $L'$ as an Euler product  provides an analytic continuation for $\Re (s) > 1/2 + \epsilon$.    The product is convergent and never zero
because its logarithm is finite,   thus there are no zeros to the right
of the critical line since $\epsilon$ can be taken arbitrarily small.   
\end{proof}

\end{theorem}

\begin{corollary}
\label{cor1}
The Dirichlet $L$-function built on the actual primes $\Primes$ is known to satisfy a functional equation that relates $L(s,\chi)$ to $L(1-s, \bar{\chi})$.    
Thus  Theorem \ref{RH1}  implies that the Generalized Riemann Hypothesis for non-principal characters is true with probability equal to $1$.
\end{corollary}

\begin{remark}
Define $\CP_{gs} $  as  the ensemble of states $\Primesp$  where $p'_n$  is a random integer satisfying 
\beq
\label{GS1}
p_n \leq p'_n \leq p_n + K,   ~~~~~  p'_n = p_n ~ {\rm mod} ~ k
\eeq
where 
$K$ is an integer.        Grosswald and Schnitzer proved that the functions $L' (s, \chi )$  can be analytically continued 
to $\Re (s) > 0$  and  remarkably have the {\it same zeros} as the $L$-function  \eqref{EPF} inside the critical strip \cite{Grosswald}.       Corollary \ref{cor1}  implies 
that all these random  $L'$-functions  based on $\CP_{gs}$   satisfy the Riemann Hypothesis with probability equal to $1$ if  $\CP_{gs} \subset  \CP$.
\end{remark}

\def\Ei{{\rm Ei}}

\section{Principal Case}

We now consider the case of the principal character of modulus $k$,  where by definition 
$\chi (n) =1$ if $(n,k)=1$,  otherwise $\chi (n) =0$.      As explained in the Introduction,  we are interested in the series
\beq
\label{BNprime}
B'_N  (t) =   \sum_{{n=1} \atop {(p'_n , k)=1}}^N \cos (t \log p'_n ) 
\eeq
where $t\neq 0$.     
As in the previous section,  this is can be modeled as 
\beq
\label{BNpz}
B'_N (t) =  \sum_{n=3}^{p'_N}   \znk  \cos (t \log n)   
\eeq

\begin{theorem}  
\label{CLTB}    
 In the limit of large $t\to \infty$,      
\beq 
\label{CLTforB}
    \sqrt{\frac{1 +  \tfrac{\log \log N}{\log N} }{s^2 N}}  ~\Bigl(  B'_N (t)  - m_N (t) \Bigr)    ~~ \dist ~~ \CN (0,1) 
\eeq
where $s^2 =  \varphi (k)/2k$  and 
\beq
\label{mNexact}
m_N (t) \approx  \Re \Bigl( \Ei\Bigl(  (1+it)\log (N \log N) \Bigr) \Bigr)
\eeq

\begin{proof}  As in the non-principal case of the last section,  the proof is based on the Lyapunov CLT.   
Let $\mu_n, \sigma_n$  be the mean and standard deviation of each term in the sum \eqref{BNpz}.    
Then
\beq
\label{mNB}
m_N (t)  = \sum_n  \mu_n =   \sum_{ {n=3} \atop {(n,k) =1}}^{p'_N}    \inv{\log n}  \cos (t \log n)   \approx  \frac{\varphi (k)}{k}   \int_3^{p'_N}    \frac{du}{\log u}   \cos (t \log u)  
\eeq
The above integral can be expressed in terms of the exponential integral function $\Ei$:  
\beq
\label{Eiint}
\int^x  \frac{du}{\log u}  \cos (t \log u)   =    \Re \(  \Ei\[  (1+it) \log x \] \)
\eeq
Using $p'_N \approx N \log N$,     we obtain  \eqref{mNexact}.  

Next let us turn to $s_N^2$:   
\barray
\nonumber
s_N^2 (t)  &=&  \sum_{n=3}^{p'_N}  \sigma_n^2   =  
\sum_{n=3}^{p'_n}  \Ex [\znk^2 \cos^2 (t \log n)]  -  \mu_n^2 
\\ \label{sNB}  
&=&  \sum_{ {n=3} \atop {(n,k) =1}}^{p'_N}  \( \inv{\log n}   -   \inv{\log^2 n}  \)  \cos^2 ( t \log n)   
\earray
We can neglect the $1/\log^2 n$  term since it is of lower order.    Approximating the sum by an integral as in \eqref{mNB},  one has   
\beq
\label{sNLi}   
s_N^2 (t)   \approx   \frac{\varphi (k)}  {2k}   \Bigl(  {\rm Li } (p'_N )   +  \Re \( \Ei [(1+2it) \log p'_N ] \)   \Bigr)
\eeq
The factor of $1/2$ in the leading ${\rm Li}$ term is a reflection that the average of $\cos^2$ is $1/2$.    
In the limit of large $t$ the $\Ei$ term can be neglected since it is  smaller by a factor of $O(1/t)$  (see the approximation in \eqref{mNapprox}).    Again using $p'_N \approx N \log N$,    
\beq
\label{sNlarget}
\lim_{t \to \infty}   s_N^2 (t)   \approx    
\frac{\varphi (k)}{2 k}  \(   \frac{N}{1 +  \tfrac{\log\log N}{\log N}}  \) 
\eeq

For the same reasons as in Theorem \ref{CLTC},   the Lyapunov condition \eqref{LyCond}   is satisfied.  
The theorem then follows from  the CLT \eqref{Lyapunov}.    
\end{proof}
\end{theorem}  

Note that for fixed $N$,   $\lim_{t\to \infty} m_N (t) =0$  (see the approximation \eqref{mNapprox} below).  
In Figure \ref{CLTBfig}   we provide numerical evidence for  Theorem \ref{CLTB}.   
As for the non-principal case,   for the  state  $\Primesp$  corresponding to  the actual primes $\Primes$, 
the series is well within the  bell curve,  namely  the LHS of \eqref{CLTforB}    for the original series $B_N (t)$ equals  $-0.280$.

\begin{remark}
Theorem \ref{CLTB}   is similar,  but not identical,  to a theorem of Kac.    For the latter,   the $p'_n$  are the true primes
and thus not random.    Rather, randomness is introduced by making $t$ a random variable,   in contrast to Theorem \ref{CLTB} 
where $t$ is not random and fixed.   Kac'  CLT is valid for $t\in [T, 2T]$  in the limit $T\to \infty$.   
\end{remark}  

\begin{figure}[t]
\centering\includegraphics[width=.5\textwidth]{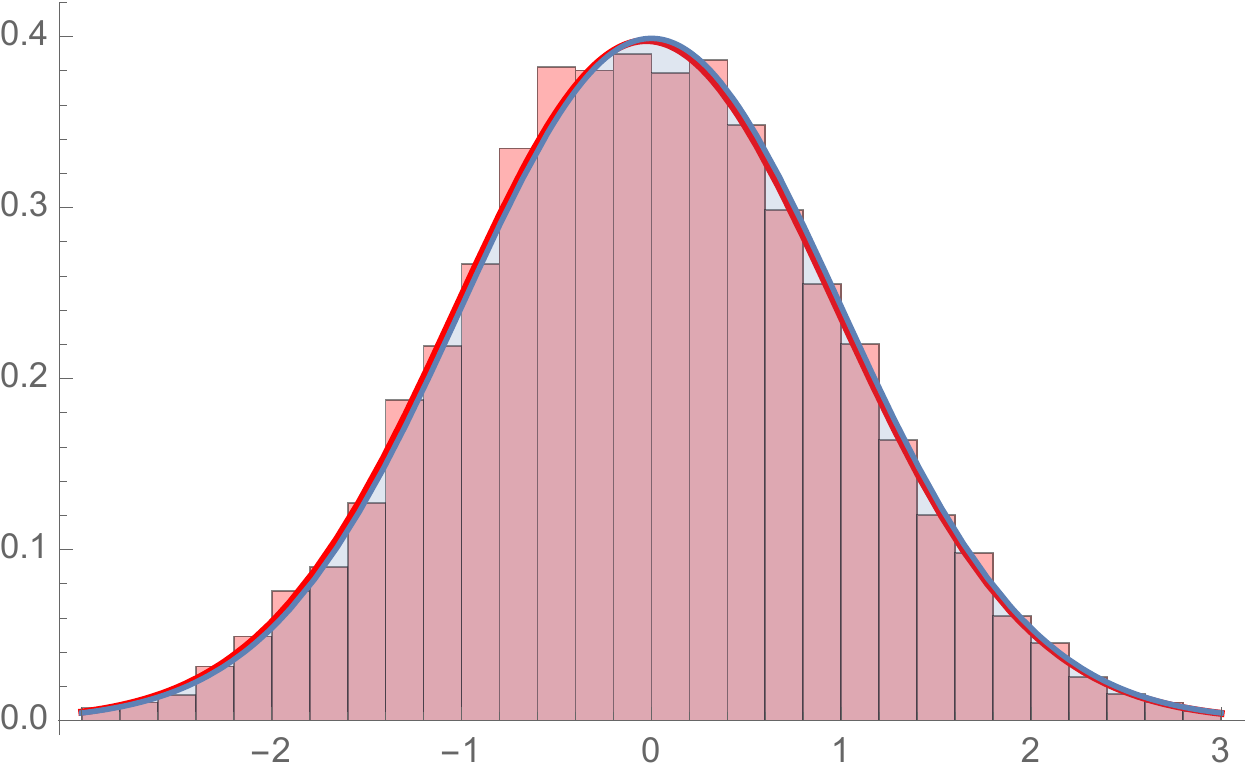}
\caption{ 
Numerical evidence for Theorem \ref{CLTB}.     We fixed $N=5,000$  and $t=1000.$
Displayed is a normalized histogram for $10,000$ states $\Primesp$.     The red curve is the  fit to the data,  which is the normal 
distribution $\CN(-0.02681, 1.00325)$.   The  blue curve is the prediction $\CN(0,1)$. }
\label{CLTBfig}
\end{figure}

\begin{theorem}   If  $t>\sqrt{N}$,    then with probability equal to one,  
\beq
\label{BsqrtN}
B'_N (t)  =  O(N^{1/2+ \epsilon})
\eeq
for any $\epsilon >0$ 
in the limit of large $N$.  
\begin{proof} 
For large $t$  and $N$,  to a very good approximation
\beq
\label{mNapprox}
m_N (t)  \approx   \frac{\varphi (k)}{k}  \(   \frac{N}{1 +  \tfrac{\log\log N}{\log N}}  \)  \( \frac{t}{1+t^2} \)   \sin \( t \log(N\log N) \)  
\eeq
If $t > \sqrt{N}$,  then  $m_N (t)  =  O(\sqrt{N})$.   
Using this,  and repeating the arguments of Theorem \ref{CNsqrt}  proves the theorem.   

\end{proof}
\end{theorem}

\begin{remark}
The above theorem implies that the Riemann Hypothesis is true with probability equal to one in the limit
$t\to \infty$.      The argument is the same as in Theorem \ref{RH1}.      
The condition $t \to \infty$ makes this a weaker statement than in the non-principal case.     
In order to deal with finite $t$,  it was proposed in \cite{Gonek,EPF1,EPFzeta}  that  a truncated Euler product 
is a good approximation to the $\zeta$ function and can be used to study the Riemann Hypothesis.   
Namely the following formula is valid:
\beq
\label{zetaRH}
\zeta (s) = \prod_{n=1}^{N(t) }  \(  1 - \inv{p_n^s} \)^{-1}   \exp \( R_{N} (s) \) 
\eeq
where $N(t) \sim  t^2$  and  $R_N (s) \sim  1/t^{2\sigma -1}$.     Thus  in the limit $t\to \infty$,  $R_N$ can be neglected 
if $\sigma  >1/2$.      The above formula would rule out zeros to the right of the critical line since the RHS is then never zero.

\end{remark}

\section*{Acknowledgments}

I   wish to thank Guilherme Fran\c ca,  Steve Gonek,  and Nicol\'as Morales-Dur\'an  for discussions.

\vfill\eject


\begin{thebibliography}{99}

\bibitem{EPF1}    G.  Fran\c ca and A.  LeClair,  
``On the validity of the Euler Product inside the critical strip",   arXiv:1410.3520 [math.NT].  


\bibitem{EPFchi}   G.  Fran\c ca and A.  LeClair, 
``Some Riemann  Hypotheses from Random Walks over Primes", 
arXiv:1509.03643 [math.NT]

\bibitem{EPFzeta}    A.  LeClair
``Riemann Hypothesis and Random Walks:  the Zeta case", 
arXiv:1601.00914 [math.NT].   

\bibitem{Grosswald}
E.  Grosswald and F. J. Schnitzer,  
{\it A class of modified $\zeta$ and $L$-functions}, 
Pacific. Jour.  Math. {\bf 74} (1978) 357.  

\bibitem{Cramer} 
H.  Cram\'er, 
``On the order of magnitude of the difference between consecutive prime numbers", 
Acta. Arith. {\bf 2}  (1936)  23.  

\bibitem{Kac}  M.  Kac,
{\it Statistical Independence in Probability,  Analysis and Number Theory,}
The Mathematical Association of America,  New Jersey,  1959.   


\bibitem{Gonek}    S.  M.  Gonek,  
``Finite Euler products and the Riemann Hypothesis",  
Trans.  Amer.  Math.  Soc.  {\bf 364} (2011)  2157.


\end{thebibliography}
\end{document}